\documentclass{article}
\usepackage{amsfonts}
\usepackage{graphics}

\usepackage{amsmath}

\newenvironment{proof}[1][Proof]{\noindent\textbf{#1 :} }{\ \rule{0.4em}{0.4em}}

\begin{document}

\title{Localizing Volatilities\footnote{I thank Marc Yor for providing key ideas for the
elaboration of this work and for all his precious comments. I also
thank H\'{e}lyette Geman for all her helpful and useful remarks.}}
\author{Marc Atlan \\
Laboratoire de Probabilit\'{e}s\\
Universit\'{e} Pierre et Marie Curie}
\date{First Draft: March 9 2004; Last Draft: April 11 2006}
 \maketitle

\begin{abstract}

We propose two main applications of Gy\"{o}ngy (1986)'s construction
of inhomogeneous Markovian stochastic differential equations that
mimick the one-dimensional marginals of continuous It\^{o}
processes. Firstly, we prove Dupire (1994) and Derman and Kani
(1994)'s result. We then present Bessel-based stochastic volatility
models in which this relation is used to compute analytical formulas
for the local volatility. Secondly, we use these mimicking
techniques to extend the well-known local volatility results to a
stochastic interest rates framework.

\end{abstract}

\section{Introduction}

It has been widely accepted for at least a decade that the option
pricing theory of Black and Scholes (1973) and Merton (1973) has
been inconsistent with option prices. Actually, the model implies
that the informational content of the option surface is one
dimensional which means that one could construct the prices of
options at all strikes and maturities from the price of any single
option. It has also been shown that unconditional returns show
excess kurtosis and skewness which are inconsistent with normality.
Special attention was given to implied volatility smile or skew, but
research has concentrated on implied Black and Scholes volatility
since it has become the unique way to price vanilla options.
Accordingly, option prices are often quoted by their implied
volatility. Nevertheless, this method is unsuitable for more
complicated exotic options and options with early exercise features.
To explain in a self-consistent way why options with different
strikes and maturities have different implied volatilities or what
one calls the
volatility smile, one could use stochastic volatility models (eg. Heston (1993) or Hull and White (1987))\\
\indent Given the computational complexity of stochastic volatility
models and the difficulty of fitting their parameters to the market
prices of vanilla options, practitioners found a simpler way to
price exotic options consistently with the volatility smile by using
local volatility models as introduced by Dupire (1994) and Derman
and Kani (1994). Local volatility models have the advantage to fit
the implied volatility surface; hence, when pricing an exotic
option, one feels comfortable hedging through the stock and vanilla options markets.\\
\indent In the last twenty years, academics and practitioners have
been primarily interested in building models that describe well the
behavior of an asset whether it is equity, FX, Credit, Fixed-Income
or Commodities and very rarely models that specify any cross-asset
dependency. For all cross-asset derivative products, this dependency
modifies the model one should use or at least the calibration
procedure. Certainly, models that incorporate a dependency on other
asset classes than a specific underlying need to be recalibrated as
soon as the other asset classes become random, in particular in the
fast growing hybrid industry where it is necessary to model several assets.\\
 \indent The remainder of the paper is organized as follows. Section
 2 recalls preliminary results on Bessel processes and states mimicking properties of
continuous It\^{o} processes exhibited by Gy\"{o}ngy (1986) and
Krylov (1985). Section 3 recalls well-known results of Dupire (1994)
and Derman and Kani (1994) on local volatility, gives a proof of the
existence of a local volatility model that mimicks a stochastic
volatility one based on Gy\"{o}ngy (1986) theorem. Section 4 gives
examples of stochastic volatility models where a local volatility
can be computed. Those examples are based on remarkable properties
of Bessel processes such as scaling properties. In order to extend
the class of volatility models (where closed-form formulas can be
obtained), we propose a general framework in which the volatility
diffusion is a general deterministic time and space transformation
of Bessel processes. Analytical computations are proposed in cases
where the volatility diffusion is independent from the stock price
diffusion as well as in cases where they are correlated. Section 5
applies the results of Section 3 to the case of stochastic interest
rates and more generally shows how Gy\"{o}ngy (1986) theorem can be
applied to construct a local volatility model in a deterministic
interest rate framework, starting from a stochastic volatility model
with stochastic rates. Finally, Section 6 concludes our work and
presents an important open question on mimicking the laws of It\^{o}
processes.

\section{Preliminary Mathematical Results}
\subsection{Bessel and CIR Processes}

Let $(R_t,t\geq 0)$ denote a Bessel process with dimension $\delta$,
starting from 0 and $(\beta_t,t\geq 0)$ an independent brownian
motion from $(B_t,t\geq 0)$ the driving brownian motion. Let us
recall that $R^2_t$ solves the following SDE:
\begin{equation*}
dR^2_t=2R_t dB_t+\delta dt
\end{equation*}
and let us now define :
\begin{eqnarray*}
I_t=\int_0^t R_s d\beta_s\indent and\indent A_t=\int_0^t R^2_s ds
\end{eqnarray*}

Then, the one-dimensional marginals of $(A_t, I_t)$ are at least in
theory well-identified, via Fourier-Laplace expressions, and are
closely related with the so-called L\'{e}vy area formula (see
L\'{e}vy (1950), Williams (1976), Gaveau (1977), Yor (1980), Chapter
2 of Yor (1992) and many other references). Here we simply recall,
for our purposes
the formulae:\\
$\forall(\alpha,\beta)\in\mathbb{R}^2$
\begin{equation}
\mathbb{E}\bigg[\exp\big(i\alpha
I_t-\frac{\beta^2}{2}A_t\big)\bigg]=\big(\cosh(t\sqrt{\alpha^2+\beta^2})\big)^{-\frac{\delta}{2}}
{\label{trans1}}
\end{equation}
as well as:\\ $\forall(a,b)\in\mathbb{R}_+\times\mathbb{R}$
\begin{equation}
\mathbb{E}\bigg[\exp\big(-a
R^2_t-\frac{b^2}{2}A_t\big)\bigg]=\big(\cosh(b
t)+\frac{2a}{b}\sinh(b t)\big)^{-\frac{\delta}{2}} {\label{trans2}}
\end{equation}
a formula that we shall use later. Some developments for the law
of $A_t$ are given, e.g. in Pitman and Yor (2003).\\
For a Bessel process of dimension $\delta$ starting at $x$, one gets
the following formula:
\\ $\forall(a,b)\in\mathbb{R}_+\times\mathbb{R}$
\begin{eqnarray*}
\mathbb{E}_x\bigg[\exp\big(-a
R^2_t-\frac{b^2}{2}A_t\big)\bigg]=\big(\cosh(b
t)+\frac{2a}{b}\sinh(b t)\big)^{-\frac{\delta}{2}}\times
\\\exp\bigg(-\frac{x^2b}{2}\frac{\sinh(bt)+\frac{2a}{b}\cosh(bt)}{\cosh(bt)+\frac{2a}{b}\sinh(bt)}\bigg)
\end{eqnarray*}

Let us now present a scaling property of the Bessel process with
respect to conditioning, which is important in the sequel.

\newtheorem{green}{Proposition}[section]
\begin{green}
{\label{prop1}} For any Bessel process $R_t$ with dimension
$\delta$, we have:
\begin{equation}
\mathbb{E}\bigg[R^2_t | \int_{0}^{t}R^2_s
ds\bigg]=\frac{2}{t}\int_{0}^{t}R^2_s ds {\label{condun}}
\end{equation}
\end{green}

\newtheorem{cond1}[green]{Remark}
\begin{cond1}
This result is in fact a very particular case of a more general
result involving only the scaling property of the process
$(R^2_t,t\geq0)$, see, e.g, Pitman and Yor (2003). But, for the sake
of completeness, we shall give a direct proof of (\ref{condun})
below:
\end{cond1}

\begin{proof}
From the scaling property of $(R^2_t,t\geq0)$, we deduce that for
every $f\in\mathcal{C}^1(\mathbb{R},\mathbb{R}_{+})$, with bounded
derivative, we have:
\begin{equation*}
\mathbb{E}\bigg[f\bigg(\int_{0}^{t}R^2_s
ds\bigg)\bigg]=\mathbb{E}\bigg[f\bigg(t^2\int_{0}^{1}R^2_s
ds\bigg)\bigg]
\end{equation*}
We then differentiate both sides with respect to $t$, to obtain:
\begin{eqnarray*}
\mathbb{E}\bigg[f'\bigg(\int_{0}^{t}R^2_s
ds\bigg)R^2_t\bigg]&=&\mathbb{E}\bigg[f'\bigg(t^2\int_{0}^{1}R^2_s
ds\bigg)(2t)\int_{0}^{1}R^2_s ds\bigg]\\
&=&\mathbb{E}\bigg[f'\bigg(\int_{0}^{t}R^2_s
ds\bigg)\frac{2}{t}\int_{0}^{t}R^2_s ds\bigg]
\end{eqnarray*}
Since this identity is true for every bounded Borel function $f'$,
the identity (\ref{condun}) follows.
\end{proof}

\newtheorem{cond1r2}[green]{Remark}
\begin{cond1r2}
We now check that formula (\ref{condun}) can be obtained directly as
a consequence of formula (\ref{trans2}): differentiating
(\ref{trans2}) both sides with respect to a and taking $a=0$, we
obtain:
\begin{equation*}
\mathbb{E}\bigg[R^2_t\exp\big(-\frac{b^2}{2}A_t\big)\bigg]=\frac{\delta}{(\cosh(bt))^{\frac{\delta}{2}+1}}\big(\frac{1}{b}\sinh(bt)\big)
\end{equation*}
while, taking $a=0$ in (\ref{trans2}), and differentiating both
sides with respect to b, we get:
\begin{equation*}
b\mathbb{E}\bigg[A_t\exp\big(-\frac{b^2}{2}A_t\big)\bigg]=\frac{\delta
t}{2(\cosh(bt))^{\frac{\delta}{2}+1}}\sinh(bt)
\end{equation*}
and the identity (\ref{trans2}) follows from the comparison of these
last two equations.
\end{cond1r2}

A reason why squared Bessel processes play an important role in
financial mathematics is that they are connected to models used in
finance. One of these models is the Cox, Ingersoll and Ross (1985)
CIR family of diffusions which are solutions of the following kind
of SDEs:
\begin{equation}
dX_t=(a-bX_t)dt+\sigma\sqrt{|X_t|}dW_t {\label{cir}}
\end{equation}
with $X_0=x_0>0$, $a\in\mathbb{R}_+$, $b\in\mathbb{R}$, $\sigma>0$
and $W_t$ a standard brownian motion. This equation admits a unique
strong (that is to say adapted to the natural filtration of $W_t$)
solution that takes values in $\mathbb{R}_+$.

One is now interested in the representation of a CIR process in
terms of a time-space transformation of a Bessel process:
\newtheorem{BCIR2}[green]{Lemma}
\begin{BCIR2}
A CIR Process $X_t$ solution of equation (\ref{cir}) can be
represented in the following form:
\begin{equation}
X_t=e^{-bt} R^2_{\frac{\sigma^2}{4b}(e^{bt}-1)}
\end{equation}
where $R$ denotes a Bessel process starting from $x_0$ at time $t=0$
of dimension $\delta=\frac{4a}{\sigma^2}$
\end{BCIR2}

\begin{proof}
This lemma results from simple properties of squared Bessel
processes that can be found in Revuz and Yor (2001), Pitman and Yor
(1980, 1982).
\end{proof}\\

This relation has been widely used in finance, for instance in Geman
and Yor (1993) or Delbaen and Shirakawa (1996).


\subsection{Mimicking Theorems}

A common topic of interest of Krylov and Gy\"{o}ngy respectively in
Krylov (1985) and Gy\"{o}ngy (1986) is the construction of
stochastic differential equations whose solutions mimick certain
features of the solutions of It\^{o} processes. The construction of
Markov martingales that have specified marginals was studied by
Madan and Yor (2002). Bibby, Skovgaard and S{\o}rensen (2005) as
well as Bibby and S{\o}rensen (1995) proposed construction of
diffusion-type models with given marginals.

Let us now consider an It\^{o} differential equation of the form:
\begin{equation}
\xi_t=\int_{0}^{t}\delta_s dW_s +\int_{0}^{t}\beta_s ds
\end{equation}
where $W_t$ is a $\mathcal{F}_t$-Brownian motion of dimension k,
$(\delta_t)_{t\in\mathbb{R}_+}$ and $(\beta_t)_{t\in\mathbb{R}_+}$
are bounded $\mathcal{F}_t$-adapted processes that belong
respectively to $\mathbf{M}_{n,k}(\mathbb{R})$, the space of
$n\times k$ real matrices and to $\mathbb{R}^n$.

\newtheorem{green2}[green]{Definition}
\begin{green2}[Green Measure]
Considering two stochastic processes $X_t$, valued in $\mathbb{R}^n$
and $\gamma_t$, with $\gamma_t>0$, one defines the Green measure
$\mu_{X,\gamma}$ by:
\begin{equation}
\mu_{X,\gamma}(\Gamma)=\mathbb{E}\bigg[\int_{0}^{\infty}\mathbf{1}_{\Gamma}(X_t)\exp\big(-\int_{0}^{t}\gamma_s
ds\big) dt\bigg]
\end{equation}
where $\Gamma$ is any borel set of $\mathbb{R}^n$
\end{green2}

\newtheorem{greenr}[green]{Remark}
\begin{greenr}
The stochastic process $\gamma_t$ is called the killing rate
\end{greenr}

\newtheorem{kryl}[green]{Theorem}
\begin{kryl}[Krylov]

If $\xi_t$ is an It\^{o} process defined as previously and
satisfying the uniform ellipticity condition: $\exists
\lambda\in\mathbb{R}_{+}^{*}$ such as $\delta\delta^{*}\geq\lambda
I_n$\\ as well as the lower boundedness condition:
$\exists\alpha\in\mathbb{R}_{+}$ such as $\gamma>\alpha$,\\ then
there exist deterministic functions
$\sigma:\mathbb{R}^n\rightarrow\mathbf{M}_{n}(\mathbb{R})$ ,
$b:\mathbb{R}^n\rightarrow\mathbb{R}$ and
$c:\mathbb{R}^n\rightarrow\mathbb{R}_+$ such that the following SDE:
\begin{eqnarray*}
dx_t&=&\sigma(x_t)dW_t+b(x_t)dt\\x_0&=&0
\end{eqnarray*}
has a weak solution satisfying:
\begin{eqnarray*}
\forall\Gamma\in\mathcal{B}(\mathbb{R}^n)\\
\mathbb{E}\bigg[\int_{0}^{\infty}\mathbf{1}_{\Gamma}(\xi_t)\exp\big(-\int_{0}^{t}\gamma_s
ds\big)
dt\bigg]&=&\mathbb{E}\bigg[\int_{0}^{\infty}\mathbf{1}_{\Gamma}(x_t)\exp\big(-\int_{0}^{t}c(x_s)
ds\big) dt\bigg]
\\ie: \mu_{\xi,\gamma}(\Gamma)&=&\mu_{X,c}(\Gamma)
\end{eqnarray*}
\end{kryl}

\begin{proof}
See Krylov (1985)
\end{proof}

\newtheorem{SDE}[green]{Definition}
\begin{SDE}[Weak Solution]
The stochastic differential equation
\begin{eqnarray}
dX_t&=&f(t,X_t)dW_t+g(t,X_t)dt {\label{sde1}}
\\X_0&=&0 {\label{sde2}}
\end{eqnarray}
is said to have a weak solution if there exist a probability space
$(\Omega, \mathcal{F},\mathbb{P})$ and an $\mathcal{F}_t$-Brownian
motion with respect to which there exists  a $\mathcal{F}_t$-adapted
stochastic process $\overline{X}_t$ that satisfies (\ref{sde1}) and
(\ref{sde2}).
\end{SDE}

A natural question asked and answered by Gy\"{o}ngy  is whether it
is possible to find the solution of an SDE with the same
one-dimensional marginal distributions as an It\^{o} process. The
answer is stated below:
\newtheorem{gyon}[green]{Theorem}
\begin{gyon}[Gy\"{o}ngy]
If $\xi_t$ is an It\^{o} process satisfying the uniform ellipticity
condition: $\exists \lambda\in\mathbb{R}_{+}^{*}$ such as
$\delta\delta^{*}\geq\lambda I_n$\\ then there exist bounded
measurable functions
$\sigma:\mathbb{R}_{+}\times\mathbb{R}^n\rightarrow\mathbf{M}_{n,n}(\mathbb{R})$
and $b:\mathbb{R}_{+}\times\mathbb{R}^n\rightarrow\mathbb{R}$
defined by:
\begin{eqnarray*}
\forall(t,x)\in\mathbb{R}_{+}\times\mathbb{R}^n\\
\sigma(t,x)&=&\bigg(\mathbb{E}\big[\delta_t\delta^{*}_t|\xi_t=x\big]\bigg)^{\frac{1}{2}}\\
b(t,x)&=&\mathbb{E}\big[\beta_t|\xi_t=x\big]
\end{eqnarray*}
\\such that the following SDE:
\begin{eqnarray*}
dx_t&=&\sigma(t,x_t)dW_t+b(t,x_t)dt\\x_0&=&0
\end{eqnarray*}
has a weak solution with the same one-dimensional marginals as
$\xi$.
\end{gyon}

\begin{proof}
See Gy\"{o}ngy (1986)
\end{proof}

\newtheorem{gykry}[green]{Remark}
\begin{gykry}
Two kinds of mimicking features of a general It\^{o} process were
illustrated in this section . With Krylov, we were able to construct
a Markov homogeneous process solution of an SDE, that has the same
Green measure than the It\^{o} process. Using Gy\"{o}ngy's results,
we were able to build a time-inhomogeneous Markov process solution
of an SDE that has the same one-dimensional marginals as the It\^{o}
process.
\end{gykry}

A possible extension to the above mimicking property is to consider
a real It\^{o} process $\xi$ driven by a multidimensional Brownian
motion and obtain a new mimicking result useful for the remainder of
the paper; the proof is straightforward from Gy\"{o}ngy (1986)
proof. Let $\xi$ be as follows :
\begin{equation}
\xi_t=\int_{0}^{t}<\delta_s,dW_s>+\int_{0}^{t}\beta_s ds
{\label{xi}}
\end{equation}
where $W_t$ is a $\mathcal{F}_t$-Brownian motion of dimension k,
$(\delta_t)_{t\in\mathbb{R}_+}$ and $(\beta_t)_{t\in\mathbb{R}_+}$
are bounded $\mathcal{F}_t$-adapted processes that belong
respectively to $\mathbb{R}^k$ and to $\mathbb{R}$.

\newtheorem{gyon2}[green]{Theorem}
\begin{gyon2}
If $\xi_t$ is an It\^{o} process defined as in (\ref{xi}) satisfying
the uniform ellipticity condition: $\exists
\lambda\in\mathbb{R}_{+}^{*}$ such as
$\delta\delta^{*}\geq\lambda$\\ then there exist bounded measurable
functions
$\sigma:\mathbb{R}_{+}\times\mathbb{R}\rightarrow\mathbb{R}$ and
$b:\mathbb{R}_{+}\times\mathbb{R}\rightarrow\mathbb{R}$ defined by:
\begin{eqnarray*}
\forall(t,x)\in\mathbb{R}_{+}\times\mathbb{R}\\
\sigma(t,x)&=&\bigg(\mathbb{E}\big[\delta_t\delta^{*}_t|\xi_t=x\big]\bigg)^{\frac{1}{2}}\\
b(t,x)&=&\mathbb{E}\big[\beta_t|\xi_t=x\big]
\end{eqnarray*}
\\such that the following SDE:
\begin{eqnarray*}
dx_t&=&\sigma(t,x_t)dW_t+b(t,x_t)dt\\x_0&=&0
\end{eqnarray*}
has a weak solution with the same one-dimensional marginals as
$\xi$.
\end{gyon2}

\section{Generalities on Local Volatility}
\subsection{Fokker-Planck Equation}

Let us assume that:
\begin{equation}
\frac{dS_t}{S_t}=r(t) dt + \sigma(t,S_t) dW_t {\label{bsloc}}
\end{equation}
where $r$ and $\sigma$ are deterministic functions, $\sigma$ is
usually called the local volatility. Under the local volatility
dynamics, option prices satisfy the following PDE:
\begin{equation}
\frac{\partial V}{\partial
t}+\frac{{\sigma^2(t,S_t)}}{2}S^2\frac{\partial^2 V}{\partial
S^2}+r(t)S\frac{\partial V}{\partial S}-r(t)V=0 {\label{PDE}}
\end{equation}
and terminal condition $V(S,T)=C(S,T)=PayOff_T(S)$.\\
If we consider call options, we would get $V(S,T)=(S-K)_{+}$. It has
been proved that one can obtain a forward PDE for $C(K,T)$ instead
of fixing $(K,T)$ and obtaining a backward PDE for $C(S,t)$. To get
the Forward PDE equation, one could just differentiate (\ref{PDE})
twice with respect to the strike K and then get the same PDE, with
variable $\phi=\frac{\partial^2 C}{\partial K^2}$ and terminal
condition $\delta(S-K)$. $\phi$ is the transition density of S and
is also the Green function of (\ref{PDE}). It follows that $\phi$ as
a function of $(K,T)$ satisfies the Fokker-Planck PDE:
\begin{equation*}
\frac{\partial\phi}{\partial T}-\frac{\partial^2}{\partial
K^2}\big(\frac{{\sigma^2(T,K)}}{2}K^2\phi\big)+r(T)\frac{\partial
}{\partial K}(K\phi)+r(T)K=0
\end{equation*}
Now, integrate twice this equation taking into account the boundary
conditions, one obtains the Forward Parabolic PDE equation:
\begin{equation}
\frac{\partial C}{\partial
T}-\frac{{\sigma^2(T,K)}}{2}K^2\frac{\partial^2 C}{\partial
K^2}+r(T)K\frac{\partial C}{\partial K}=0 {\label{fPDE}}
\end{equation}
with initial condition $C(K,0)=(S_0-K)_{+}$. Hence, one obtains
Dupire (1994) equation
\begin{equation}
\sigma^2(T,K)=\frac{\frac{\partial C}{\partial
T}+r(T)K\frac{\partial C}{\partial
K}}{\frac{1}{2}K^2\frac{\partial^2 C}{\partial K^2}}
{\label{dupire}}
\end{equation}
\\Moreover, if one expresses
the option price as a function of the forward price, one would write
a simpler expression:
\begin{equation*}
\sigma^2(T,K,S_0)=\frac{\frac{\partial C}{\partial
T}}{\frac{1}{2}K^2\frac{\partial^2 C}{\partial K^2}}
\end{equation*}
where C is now a function of $(F_T,K,T)$ with
$F_T=S_0\exp\big(\int_{0}^{T}r(s) ds\big)$.

\subsection{Matching Local and Stochastic Volatilities}

A stock price diffusion with a stochastic volatility is one of the
following form:
\begin{equation}
\frac{dS_t}{S_t}=r(t) dt +\sqrt{V_t} dW_t {\label{bssto}}
\end{equation}
where $V_t$ is a stochastic process, solution of an SDE and $r(t)$ is a deterministic function of time.
(We do not yet discuss the dependence of the stock price and volatility processes, also called Leverage effect)\\
One can find a relation between the local volatility and a
stochastic volatility. First, one applies Tanaka's formula to the
stock price process:
\begin{eqnarray*}
e^{-\int_0^t r(s) ds}(S_t-x)_{+}&=&(S_0-x)_{+} -\int_0^t r(u)
e^{-\int_0^u r(s) ds}(S_u-x)_{+}du \\&&+ \int_0^t e^{-\int_0^u r(s)
ds}\mathbf{1}_{\{S_u>x\}}dS_u+\frac{1}{2}\int_0^t e^{-\int_0^u r(s)
ds} dL_u^x(S)
\end{eqnarray*}
Assuming that $(e^{-\int_0^t r(s) ds}S_t, t\geq 0)$ is a true
martingale, then\\ $(\int_0^t\mathbf{1}_{\{S_u>x\}}d(e^{-\int_0^u
r(s) ds}S_u) , t\geq 0)$ is a martingale and one gets:
\begin{eqnarray*}
\mathbb{E}[e^{-\int_0^t r(s)
ds}(S_t-x)_{+}]&=&\mathbb{E}[(S_0-x)_{+}] +x\int_0^t \mathbb{E}[r(u)
e^{-\int_0^u r(s) ds}\mathbf{1}_{\{S_u>x\}}] du \\&&+\frac{1}{2}
\mathbb{E}\bigg[\int_0^t e^{-\int_0^u r(s) ds} dL_u^x(S)\bigg]
\end{eqnarray*}
Then, differentiating the previous relation and using Fubini
theorem, one obtains:
\begin{equation}
d_t C(t,x)= x\mathbb{E}[r(t) e^{-\int_0^t r(s)
ds}\mathbf{1}_{\{S_t>x\}}] dt+ \frac{1}{2}\mathbb{E}[ e^{-\int_0^t
r(s) ds}dL_t^a(S)]
\end{equation}
where $C(t,x)=\mathbb{E}[e^{-\int_0^t r(s) ds}(S_t-x)_{+}]$ Using a
classical characterization of the local time of any continuous
semi-martingale:
\begin{equation}
L_t^x(S)=\lim_{\epsilon\rightarrow
0}\frac{1}{\epsilon}\int_{0}^{t}\mathbf{1}_{\{x\leq S_s<
x+\epsilon\}}d<S,S>_s
\end{equation} one gets with a permutation of the differentiation and the expectation:
\begin{equation}
d_t C(t,x)=x\mathbb{E}[r(t) e^{-\int_0^t r(s)
ds}\mathbf{1}_{\{S_t>x\}}] dt+\frac{1}{2}\lim_{\epsilon\rightarrow
0}\mathbb{E}[\frac{1}{\epsilon}\mathbf{1}_{\{x\leq S_t<
x+\epsilon\}}e^{-\int_0^t r(s) ds} V_t S^2_t] dt{\label{proofvol}}
\end{equation}
as a result of $d<S,S>_t=V_t S^2_t dt$. Now, one may write using
conditional expectations and the fact that interest rates are
assumed to be deterministic, the following identity:
\begin{equation*}
\mathbb{E}[V_t
S^2_t\mathbf{1}_{\{x\leq S_t< x+\epsilon\}}]=\mathbb{E}\big[\mathbb{E}[V_t|S_t]S_t^2\mathbf{1}_{\{x\leq S_t< x+\epsilon\}}\big]\\
\end{equation*}
From there, one easily obtains:
\begin{eqnarray*}
\lim_{\epsilon\rightarrow 0}\frac{1}{\epsilon}\mathbb{E}[V_t
S^2_t\mathbf{1}_{\{x\leq S_t< x+\epsilon\}}]&=&\lim_{\epsilon\rightarrow 0}\frac{1}{\epsilon}\mathbb{E}\big[\mathbb{E}[V_t|S_t]S_t^2\mathbf{1}_{\{x\leq S_t< x+\epsilon\}}\big]\\
&=&\mathbb{E}[V_t|S_t=x] x^2 q_t(x)
\end{eqnarray*}
where $q_t(x)$ is the value of the density of $S_t$ in $x$. Since
Breeden and Litzenberger (1978), it is well known that
$\frac{\partial^2 C(t,x)}{\partial x^2}=e^{-\int_0^t r(s) ds}
q_t(x)$. It is also known that $\frac{\partial C}{\partial x} =
-\mathbb{E}[e^{-\int_0^t r(s) ds}\mathbf{1}_{\{S_t>x\}}]$ One
finally may write:
\begin{equation}
\frac{\partial C}{\partial t} +x r(t)\frac{\partial C}{\partial
x}=\mathbb{E}[V_t|S_t=x]\frac{1}{2}x^2\frac{\partial^2 C}{\partial
x^2}{\label{dupire2}}
\end{equation}

\noindent Comparing equation (\ref{dupire}) and the above equation,
one may obtain an equation that relates local and stochastic
volatility models
\begin{equation}
\sigma^2(t,x)=\mathbb{E}[V_t|S_t=x] {\label{mimic}}
\end{equation}
Hence, we have proven that if there exists a local volatility such
as the one-dimensional marginals of the stock price with the implied
diffusion are the same as the ones of the stock price with the
stochastic volatility, then the local volatility satisfies equation
(\ref{mimic}).

Another way to prove this relation is to apply Gy\"{o}ngy (1986)
result. Since the stock price dynamics with a stochastic volatility
given by equation (\ref{bssto}) and the ones with the local
volatility given by equation (\ref{bsloc}) must have the same
one-dimensional marginals, one can apply Gy\"{o}ngy Theorem:
assuming that there exists $\lambda\in\mathbb{R}_{+}^{*}$ such that
$S^2V\geq\lambda$ we get the well-known relation between the local
and the stochastic volatilities:
\begin{equation*}
\sigma(t,S_t=x) =
\bigg(\mathbb{E}\big[V_t|S_t=x\big]\bigg)^{\frac{1}{2}}
\end{equation*}
It is important to notice that Gy\"{o}ngy gives us the existence of
such a diffusion in addition to provide an explicit way to construct
it. More generally, assuming just that the volatility process is a
general continuous semi-martingale, one can also get the same
result, and a justification for the use of local volatility models.
Hence, we obtain an illustration of Gy\"{o}ngy's result in a finance
framework. Moreover, it is shown that one can get the relation
(\ref{mimic}) without using the Forward PDE equation.

As a first remark, we should notice that if we choose $V_t$ such as
$\sqrt{V_t}=\sigma(t,S_t)$, we then obtain another direct proof of
equation (\ref{dupire}).\\

As a second remark, we can prove that if $(\widetilde{S}_t =
e^{-\int_0^t r(s) ds}S_t) , t\geq 0)$ is a strict local martingale
(which is studied in Cox and Hobson (2005) who named this market
situation a bubble), then
\begin{eqnarray*}
\mathbb{E}[\int_0^t\mathbf{1}_{\{S_u>x\}}d\widetilde{S}_u]&=&\mathbb{E}[\widetilde{S}_t-S_0]-\mathbb{E}[\int_0^t\mathbf{1}_{\{S_u\leq
x\}}d\widetilde{S}_u]\\
&=&\mathbb{E}[\widetilde{S}_t-S_0]
\end{eqnarray*}
since using Madan and Yor (2006), $(\int_0^t\mathbf{1}_{\{S_u\leq
x\}}d(e^{-\int_0^u r(s) ds}S_u) , t\geq 0)$ is a square integrable
martingale. Hence, defining
\begin{equation*}
c_{\widetilde{S}}(t)=\mathbb{E}[S_0-\widetilde{S}_t]
\end{equation*}
and assuming that $c_{\widetilde{S}}$ is a continuously
differentiable function, one obtains an extension of equation
(\ref{dupire}) that is a generalization to the case of strict local
martingales. This equation writes
\begin{equation*}
\sigma^2(t,x)= \frac{\frac{\partial C}{\partial t} +x
r(t)\frac{\partial C}{\partial
x}+c'_{\widetilde{S}}(t)}{\frac{1}{2}x^2\frac{\partial^2 C}{\partial
x^2}}
\end{equation*}

\section{Applications to the Heston (1993) model and Extensions}
\subsection{The Simplest Heston Model}

The aim of this paragraph is now to compute the local volatility not
by excerpting it from the option prices (see for instance Derman and
Kani (1994)) but by applying Gy\"{o}ngy's theorem.

Among the possible choices of stochastic volatility models, we will
consider the simplest one, given by the following SDE:
\begin{equation}
\frac{dS_t}{S_t}=W_t dB_t {\label{simple}}
\end{equation}
\begin{equation*}
S_0=1
\end{equation*}
where $(W_t)$ and $(B_t)$ are two independent one-dimensional
Brownian motions starting at 0. We do not consider any drift term in
our stock diffusion as we look at the forward price dynamics that
are driftless by construction.

To make our discussion a little more general than the model
presented in equation ({\ref{simple}}), we write ({\ref{simple}}):
\begin{equation*}
\frac{dS_t}{S_t}=|W_t|sgn(W_t) dB_t
\end{equation*}
\begin{equation*}
S_0=1
\end{equation*}
Now we define $\beta_t=\int_0^t sgn(W_s)dB_s$, another Brownian
motion which is independent of $(W_t,t\geq 0)$ and consequently of
the reflecting Brownian motion $(|W_t|,t\geq 0)$. We get the
following model:
\begin{equation*}
\frac{dS_t}{S_t}=|W_t|d\beta_t,\indent S_0=1
\end{equation*}
Now this modified form leads itself naturally to the generalization:
\begin{equation}
\frac{dS_t}{S_t}=R_t d\beta_t,\indent S_0=1 {\label{gene}}
\end{equation}
where, as in subsection 2.1, $(R_t)$ denotes a Bessel process with
dimension $\delta$ starting at $0$ and $(\beta_t)$ an independent
Brownian motion.

Let us consider a Markovian martingale ($\Sigma_t,t\geq0$), which is
the unique solution of:
\begin{equation}
\frac{d\Sigma_t}{\Sigma_t}=\sigma(t,\Sigma_t)d\beta_t
{\label{volloc}}
\end{equation}
\begin{equation*}
\Sigma_0=1
\end{equation*}
for some particular diffusion coefficient
\{$\sigma(t,x),t\geq0,x\in\mathbb{R}_{+}$\} which has the same
one-dimensional marginal distributions as $(S_t,t\geq0)$ the
solution of (\ref{gene}).

We will now use proposition \ref{prop1} to find $\sigma$, the local
volatility. We follow the notation in subsection 2.1, and introduce
a useful notation:

\begin{eqnarray}
L^{(\mu)}_t &=&I_t-\mu A_t\\
&\overset{(Law)}{=}& N\sqrt{A_t}-\mu A_t {\label{mu}}
\end{eqnarray}
where $N$ is a standard gaussian variable independent of $A_t$. Next
we remark as a consequence of (\ref{mu}) that for any fixed $t\geq0$
:
\begin{equation*}
(R_t,L^{(\mu)}_t)\overset{(Law)}{=}(R_t,N\sqrt{A_t}-\mu A_t)
\end{equation*}
and
\begin{equation*}
\mathbb{E}\big[R^2_t|L^{(\mu)}_t=l\big]=\mathbb{E}\big[\mathbb{E}(R^2_t|N,A_t)|N\sqrt{A_t}-\mu
A_t=l\big]
\end{equation*}
Since $N$ is independent of $R_t$, we obtain
\begin{equation*}
\mathbb{E}\big[R^2_t|L^{(\mu)}_t=l\big]=\mathbb{E}\big[\mathbb{E}(R^2_t|A_t)|N\sqrt{A_t}-\mu
A_t=l\big]
\end{equation*}
From (\ref{condun}), we deduce:
\begin{equation}
\mathbb{E}\big[R^2_t|L^{(\mu)}_t=l\big]=\big(\frac{2}{t}\big)\mathbb{E}\big[A_t|N\sqrt{A_t}-\mu
A_t=l\big] {\label{lmu}}
\end{equation}
Now, the computation of the expression in (\ref{lmu}) is a simple
exercise, which we present in the following form:
\newtheorem{lemu}[green]{Lemma}
\begin{lemu}
Let $X>0$ be a random variable independent from a standard gaussian
variable $N$.
Denote $Y^{(\mu)}=N\sqrt{X}-\mu X$. Then:\\
$i)$ for any $f:\mathbb{R}_{+}\rightarrow\mathbb{R}_{+}$, Borel
function, the following formula holds:
\begin{equation*}
\mathbb{E}\big[f(X)|Y^{(\mu)}=z\big]=\frac{h^{(\mu)}(f;z)}{h^{(\mu)}(1;z)}
\end{equation*}
where:
$h^{(\mu)}(f;z)=\mathbb{E}\bigg[\frac{f(X)}{\sqrt{X}}\exp\bigg(-\frac{(z+\mu
X)^2}{2X}\bigg)\bigg]$\\
$ii)$ in particular, for $f(x)=x$, one can write :
\begin{equation}
\mathbb{E}\big[X|Y^{(\mu)}=z\big]=-\Bigg(\frac{\frac{\partial
k}{\partial b}}{k}\Bigg)\big(\frac{z^2}{2},\frac{\mu^2}{2}\big)
{\label{condmu}}
\end{equation}
where
$k(a,b)=\mathbb{E}\bigg[\frac{1}{\sqrt{X}}\exp\bigg(-\big(\frac{a}{X}+bX\big)\bigg)\bigg]$
\end{lemu}

The proof of this lemma results from elementary properties of
conditioning and is left to the reader.\\ We now give a formula for
$\sigma^2(t,x)$ in terms of the law of $A_t\equiv A^{(\delta)}_t$,
by using equation (\ref{lmu}) and the above lemma. Indeed, it
follows from these results that:
\begin{equation*}
\mathbb{E}\big[R^2_t|\ln(S_t)=l\big]=-\frac{2}{t}\frac{\frac{\partial
k^t_{\delta}}{\partial
b}\big(\frac{l^2}{2},\frac{1}{8}\big)}{k^t_{\delta}\big(\frac{l^2}{2},\frac{1}{8}\big)}
\end{equation*}
where
$k^t_{\delta}(a,b)=\mathbb{E}\bigg[\frac{1}{\sqrt{A_t}}\exp\bigg(-\big(\frac{a}{A_t}+b
A_t\big)\bigg)\bigg]$.

Using the scaling property, we have
$k^t_{\delta}(a,b)=\frac{1}{t}k^1_{\delta}(\frac{a}{t^2},b t^2)$
which allows us to concentrate on $k^{\delta}(a,b)\equiv
k^1_{\delta}(a,b)$.\\
The following formula for the density $f_{\delta}$ of $A_{1}$ is
borrowed from Biane, Pitman and Yor (2001). Denoting
$h=\frac{\delta}{2}$, we have:
\begin{equation}
f_{\delta}(x)\equiv
f^{\sharp}_h(x)=\frac{2^h}{\Gamma(h)}\sum_{n=0}^{\infty}(-1)^n\frac{\Gamma(n+h)}{\Gamma(n+1)}\frac{(2n+h)}{\sqrt{2\pi
x^3}}\exp\bigg(-\frac{(2n+h)^2}{2x}\bigg) {\label{series}}
\end{equation}
For $\delta=2$, $A^{(2)}$, or equivalently $f_2(x)=f^{\sharp}_1(x)$
enjoys a symmetry property (also shown in Biane,
Pitman and Yor (2001)):\\
For any non-negative measurable function $g$
\begin{equation}
\mathbb{E}\bigg[g\large(\frac{4}{\pi^2
A^{(2)}}\large)\bigg]=\sqrt{\frac{2}{\pi}}\mathbb{E}\bigg[\frac{1}{\sqrt{A^{(2)}}}g(A^{(2)})\bigg],
\end{equation}
\begin{equation}
f^{\sharp}_1(x)=\large(\frac{2}{\pi
x}\large)^{\frac{3}{2}}f^{\sharp}_1\large(\frac{4}{\pi^2 x}\large)
\end{equation}
and
\begin{equation}
f^{\sharp}_1(x)=\pi \sum_{n=0}^{\infty}(-1)^n
\large(n+\frac{1}{2}\large)e^{-\large(n+\frac{1}{2}\large)^2\pi^2\frac{x}{2}}
\end{equation}
From formula (\ref{series}), one may compute with the change of
variables $a=\frac{\alpha^2}{2}, b=\frac{\beta^2}{2}$
\begin{eqnarray*}
k^{\delta}(a,b)
&\equiv&\mathbb{E}\bigg[\frac{1}{\sqrt{A^{(\delta)}}}\exp\bigg(-\frac{1}{2}\big(\frac{\alpha^2}{A^{(\delta)}}+\beta^2
A^{(\delta)}\big)\bigg)\bigg]\\
&=&\frac{2^h}{\Gamma(h)}\sum_{n=0}^{\infty}(-1)^n\frac{\Gamma(n+h)}{\Gamma(n+1)}\frac{2n+h}{\sqrt{2\pi}}\int_{0}^{\infty}\frac{dx}{x^2}e^{-\frac{1}{2}\big(\frac{\alpha^2+(2n+h)^2}{x}+\beta^2
x\big)}\\
&=&\frac{2^h}{\Gamma(h)}\sum_{n=0}^{\infty}(-1)^n\frac{\Gamma(n+h)}{\Gamma(n+1)}\frac{2n+h}{\sqrt{2\pi}}\int_{0}^{\infty}
e^{-\frac{1}{2}\big((\alpha^2+(2n+h)^2)x+\frac{\beta^2}{x}\big)}dx\indent
(\star)
\end{eqnarray*}
Also of importance for us, is the result:
\begin{equation}
\frac{\partial}{\partial
b}\large(k^{\delta}(a,b)\large)=\frac{-(2^h)}{\Gamma(h)}\sum_{n=0}^{\infty}(-1)^n\frac{\Gamma(n+h)}{\Gamma(n+1)}\frac{2n+h}{\sqrt{2\pi}}\int_{0}^{\infty}
\frac{e^{-\frac{1}{2}\big(\alpha_n^2
x+\frac{\beta^2}{x}\big)}}{x}dx {\label{derive}}
\end{equation}
where $\alpha_n=\sqrt{\alpha^2+(2n+h)^2}$.\\
Recall the integral representation for the Mc Donald functions
$K_{\nu}$:
\begin{equation}
K_{\nu}(z)\equiv
K_{-\nu}(z)=\frac{1}{2}\big(\frac{z}{2}\big)^{\nu}\int_{0}^{\infty}\frac{dt}{t^{\nu+1}}\exp-\big(t+\frac{z^2}{2t}\big){\label{kmu}}
\end{equation}
In particular, we have:
\begin{equation*}
K_{0}(z)=\frac{1}{2}\int_{0}^{\infty}\frac{dt}{t}e^{-\big(t+\frac{z^2}{2t}\big)}
\end{equation*}
As a consequence:
\begin{equation}
\int_{0}^{\infty}\frac{du}{u}e^{-\frac{1}{2}\big(\alpha^2
u+\frac{\beta^2}{u}\big)}=2K_{0}\bigg(\frac{\alpha\beta}{\sqrt{2}}\bigg)
{\label{macdo}}
\end{equation}
Now, plugging (\ref{macdo}) in (\ref{derive}), we obtain:
\begin{equation}
\frac{\partial}{\partial
b}\large(k^{\delta}(a,b)\large)=\frac{-(2^h)}{\Gamma(h)}\sum_{n=0}^{\infty}(-1)^n\frac{\Gamma(n+h)}{\Gamma(n+1)}\frac{2n+h}{\sqrt{2\pi}}2K_{0}\bigg(\frac{\alpha_n\beta}{\sqrt{2}}\bigg)
{\label{dkdel}}
\end{equation}
Likewise, we deduce from (\ref{kmu}) that:
\begin{equation*}
K_1(z)\equiv K_{-1}(z)=\frac{1}{z}\int_0^{\infty}dt
e^{-\big(t+\frac{z^2}{2t}\big)}
\end{equation*}
which implies
\begin{equation*}
\int_0^{\infty}du e^{-\frac{1}{2}\big(\alpha^2_n
u+\frac{\beta^2}{u}\big)}=\frac{\beta\sqrt{2}}{\alpha_n}
K_1\big(\frac{\alpha_n\beta}{\sqrt{2}}\big)
\end{equation*}
Hence, we get as a consequence of $(\star)$:
\begin{equation}
k^{\delta}(a,b)=\frac{2^h}{\Gamma(h)}\sum_{n=0}^{\infty}(-1)^n\frac{\Gamma(n+h)}{\Gamma(n+1)}\frac{2n+h}{\sqrt{\pi}}\frac{\beta}{\alpha_n}K_{1}\bigg(\frac{\alpha_n\beta}{\sqrt{2}}\bigg)
{\label{kdel}}
\end{equation}
Recalling that $\beta=\sqrt{2b}$ and that
$\alpha_n=\sqrt{2a+(2n+h)^2}$, we may now write (\ref{kdel}) and
(\ref{dkdel}) as:
\begin{eqnarray}
k^{\delta}(a,b)=\frac{2^h}{\Gamma(h)}\sum_{n=0}^{\infty}(-1)^n\frac{\Gamma(n+h)}{\Gamma(n+1)}\frac{2n+h}{\sqrt{\pi}}\frac{\sqrt{2b}}{\alpha_n}K_{1}(\alpha_n\sqrt{b})
\\\frac{\partial}{\partial
b}\large(k^{\delta}(a,b)\large)=\frac{-(2^h)}{\Gamma(h)}\sum_{n=0}^{\infty}(-1)^n\frac{\Gamma(n+h)}{\Gamma(n+1)}\frac{2n+h}{\sqrt{2\pi}}2K_{0}(\alpha_n\sqrt{b})
\end{eqnarray}

And finally, we obtain the following formula for the local
volatility:
\begin{equation}
\sigma^2(t,x=e^l)=-\frac{4}{t}\bigg(\frac{\frac{\partial
k^{\delta}}{\partial
b}}{k^{\delta}}\bigg)\big(\frac{l^2}{2t^2},\frac{t^2}{8}\big)
{\label{eqloc}}
\end{equation}

\subsection{Adding the Correlation}

We now assume a non-zero correlation between the volatility process
and the stock price process. This is a common fact in finance called
the Leverage Effect and translated by a negative correlation. For a
financial understanding of this effect, one can refer for instance
to Black (1976), Christie (1982) or Schwert (1989).

Let us define our new model for the stock price dynamics with a
Bessel process of dimension $\delta$ starting from 0 correlated to
the Brownian motion of the stock price process:
\begin{eqnarray}
\frac{dS_t}{S_t}=R_t dW_t {\label{stdyn}}\\
dR_t^2=2 R_t dW^{\sigma}_t+\delta t {\label{voldyn}}\\
d<W^{\sigma},W>_t=\rho dt\\
S_0=1\indent and\indent R_0=0
\end{eqnarray}

Then, there exists a Brownian motion $\beta$ independent of the
Bessel process such that $\forall t$:
\begin{equation*}
W_t=\rho W^{\sigma}_t+\sqrt{1-\rho^2}\beta_t
\end{equation*}
Using the previous formula, plugging it in (\ref{stdyn}) and then
inserting (\ref{voldyn}) in the new (\ref{stdyn}), one gets:
\begin{equation}
\frac{dS_t}{S_t}=\frac{\rho}{2}(dR_t^2-\delta
dt)+\sqrt{1-\rho^2}R_t d\beta_t
\end{equation}
Then using It\^{o} formula applied to $f(x)=ln(x)$
\begin{equation*}
d\ln(S_t)=\frac{dS_t}{S_t}-\frac{1}{2}R_t^2 dt
\end{equation*}
one obtains:
\begin{equation}
\ln(S_t)=\frac{\rho}{2}(R_t^2-\delta t)+\sqrt{1-\rho^2}\int_0^t
R_s d\beta_s -\frac{1}{2}\int_0^t R_s^2 ds {\label{lnst}}
\end{equation}

Let us consider as in subsection 5.1, $L_t^{(\mu)}=\int_0^t R_s
d\beta_s-\mu\int_0^t R^2_s ds$ (we are especially interested in the
case $\mu=\frac{1}{2\sqrt{1-\rho^2}}$). Since R and $\beta$ are
independent, we shall use the same notation as above.
 Particularly, $A_t$ and $I_t$ will refer to the quantities defined in subsection 2.1.\\
Now, equation (\ref{lnst}) can be rewritten as follows:
\begin{equation}
\ln(S_t)=\frac{\rho}{2}(R_t^2-\delta
t)+\sqrt{1-\rho^2}L_t^{(\frac{1}{2\sqrt{1-\rho^2}})}
\end{equation}
Since we wish to evaluate the local volatility
$\mathbb{E}\big[R^2_t|\ln(S_t)=l\big]$, we will try to compute more
generally the following quantity:
\begin{equation}
\mathbb{E}\big[R^2_t|m R_t^2+L_t^{(\mu)}=l \big]
\end{equation}
where $m$ is a real constant.

\newtheorem{correl}[green]{Remark}
\begin{correl}
We immediately see that if we take $m=0$, i.e $\rho=0$, we are back
to the previous paragraph setting.
\end{correl}

First, we see that equation (\ref{lmu}) is easily extended to the
case with correlation and we obtain:
\begin{equation}
\mathbb{E}\big[R^2_t|m R_t^2+L_t^{(\mu)}=l
\big]=\frac{2}{t}\mathbb{E}\big[A_t|m R_t^2+L_t^{(\mu)}=l \big]
\end{equation}
Before extending Lemma 4.1, one must recall that for any $t\geq0$:
\begin{equation*}
(R_t,A_t,L^{(\mu)}_t)\overset{(Law)}{=}(R_t,A_t,N\sqrt{A_t}-\mu A_t)
\end{equation*}
where $N$ is a standard gaussian variable independent of $R_t$ and
$A_t$. The following simple result will be helpful for the
remaining of the paper:\\

\newtheorem{lemu2}[green]{Lemma}
\begin{lemu2}
Let $X>0$ and $Z\geq 0$ independent from a standard gaussian
variable $N$.
Denote $Y^{(\mu)}=N\sqrt{X}-\mu X$. Then:\\
$i)$ for any Borel function
$f:\mathbb{R}_{+}^2\rightarrow\mathbb{R}_{+}$ , real number $m$ we
have the formula:
\begin{equation*}
\mathbb{E}\big[f(X,Z)|mZ+Y^{(\mu)}=z\big]=\frac{a^{(\mu,m)}(f;z)}{a^{(\mu,m)}(1;z)}
\end{equation*}
where:
$a^{(\mu,m)}(f;z)=\mathbb{E}\bigg[\frac{f(X,Z)}{\sqrt{X}}\exp\bigg(-\frac{(z+\mu
X-mZ)^2}{2X}\bigg)\bigg]$\\
$ii)$ in particular, for $f(x,y)=x$, we obtain:
\begin{equation}
\mathbb{E}\big[X|mZ+Y^{(\mu)}=z\big]=-\frac{1}{\mu\sqrt{2}}\Bigg(\frac{\frac{\partial
\alpha}{\partial
b}}{\alpha}\Bigg)\big(\frac{z}{\sqrt{2}},\frac{\mu}{\sqrt{2}},\frac{m}{\sqrt{2}}\big)
{\label{condmu2}}
\end{equation}
where
$\alpha(a,b,c)=\mathbb{E}\bigg[\frac{1}{\sqrt{X}}\exp\bigg(-\big(\frac{(a-cZ)^2}{X}+b^2X+bcZ\big)\bigg)\bigg]$
\end{lemu2}
The other fundamental result we now need, is the joint
density of $(R_t^2,\int_0^t R_s^2 ds)_{t\geq 0}$.\\
\newtheorem{lemu3}[green]{Theorem}
\begin{lemu3}
The joint distribution $g_t$ of $(R_t^2,\int_0^t R_s^2 ds)$ is
given by:
\begin{equation}
g_t(x,y)=\frac{1}{\sqrt{2\pi}\Gamma(\frac{\delta}{2})}\sum_{j=0}^{\infty}\frac{(-1)^j}{j!}x^{j+\frac{\delta}{2}-1}y^{-\frac{j}{2}-\frac{\delta}{4}-1}f_t^j(x,y)
\end{equation}
where $f_t^j$ is defined by
\begin{equation}
f_t^j(x,y)=\sum_{k=0}^{\infty}\frac{(j+\frac{\delta}{2})_k}{k!}e^{-\frac{1}{4y}[2(k+j+\frac{\delta}{4})t+\frac{x}{2}]^2}D_{\frac{\delta}{2}+j+1}\big(\frac{2(k+j+\frac{\delta}{4})t+\frac{x}{2}}{\sqrt{y}}\big)
\end{equation}
$D_{\nu}(\xi)$ is a parabolic cylinder function and $(\nu)_k$ the
Pochhammer's symbol defined by $(\nu)_k\equiv\nu(\nu+1). .
.(\nu+k-1)=\Gamma(\nu+k)/\Gamma(\nu)$
\end{lemu3}

\begin{proof}
See Ghomrasni (2004) who evaluates the Laplace transform of
(\ref{trans2}) in order to get the density function.
\end{proof}\\

\noindent For the definition and properties on the parabolic
cylinder functions, we refer to Gradshteyn and Ryzhik (2000).\\
\noindent Let us define $\alpha_t$ in the following form:
\begin{equation}
\alpha_t(a,b,c)=\mathbb{E}\bigg[\frac{1}{\sqrt{A_t}}\exp\bigg(-\big(\frac{(a-cR_t^2)^2}{A_t}+b^2A_t+bc
R_t^2\big)\bigg)\bigg]
\end{equation}
\noindent Unfortunately, there is no more scaling property as in the
zero-correlation case and we may not rewrite $\alpha_t$ as a
function of $t$ and $\alpha_1$. One can then compute the local
volatility $\sigma^{(\rho)}$ by noticing that in the case of
particular interest for us, the parameters are defined as follows:
\begin{equation*}
m=\frac{\rho}{2\sqrt{1-\rho^2}}\indent and\indent
z=\frac{l+\frac{\rho}{2}\delta t}{\sqrt{1-\rho^2}}\indent
and\indent\mu=\frac{1}{2\sqrt{1-\rho^2}}
\end{equation*}
We then obtain
\begin{equation}
\sigma^{(\rho)}(t,x=e^l)=-\sqrt{2(1-\rho^2)}\Bigg(\frac{\frac{\partial
\alpha_t}{\partial
b}}{\alpha_t}\Bigg)\big(\frac{l+\frac{\rho}{2}\delta
t}{\sqrt{2(1-\rho^2)}},\frac{1}{\sqrt{8(1-\rho^2)}},\frac{\rho}{\sqrt{8(1-\rho^2)}}\big)
\end{equation}

\subsection{From a Bessel Volatility process to the Heston Model}

The Heston (1993) model for representing a stochastic volatility
process is a particular case of the Cox, Ingersoll and Ross (1985)
stochastic process, of the form:
\begin{equation}
dV_t=\kappa(\theta-V_t)dt+\eta\sqrt{V_t}dW_t
\end{equation}
with initial condition $V_0=v_0$

Actually, it is possible to find out deterministic space and time
changes such as the law of the Heston SDE solution and the
Time-Space transformed Bessel Process are the same.

\newtheorem{tspace}[green]{Proposition}
\begin{tspace}
For every Heston SDE solution, there exist a Bessel process and two
deterministic functions $f$ and $g$ with $g$ increasing such as:
\begin{equation*}
V_t=f(t)\times R^2_{g(t)}
\end{equation*}
where $R$ denotes a Bessel Process of dimension
$\delta=\frac{4\kappa\theta}{\eta^2}$ starting from $\sqrt{v_0}$ at
time $t=0$ and $f$ and $g$ are defined by:
\begin{eqnarray*}
f(t)=e^{-\kappa t}\\
g(t)=\frac{\eta^2}{4\kappa}(e^{\kappa t}-1)
\end{eqnarray*}
\end{tspace}

\begin{proof}
It is just an application of Lemma 2.4.
\end{proof}\\

One may now apply the results of the previous sections using the
time and space transformations presented in the previous paragraph

\newtheorem{locheston}[green]{Proposition}
\begin{locheston}
Let us consider the following stochastic volatility model:
\begin{eqnarray*}
\frac{dS_t}{S_t}&=&\sqrt{v_t}d\beta_t,\indent S_{\{t=0\}}=S_0\\
v_t&=&\frac{\eta^2}{4}e^{2\kappa t}V_t,\indent v_0=\frac{\eta^2}{4}V_0 \\
dV_t&=&\kappa(\theta-V_t)dt+\eta\sqrt{V_t}dW_t\\
d<\beta_{.},W_{.}>_t&=&\rho dt
\end{eqnarray*}
where $\beta_t$ is a Brownian motion and $V_t$ is an Heston process
as defined above.
\\Then the local volatility $\tilde{\sigma}$ that gives us the expected mimicking properties, satisfies the
following equation:
\begin{equation}
\tilde{\sigma}(t,x)=\frac{\eta^2}{4}e^{\kappa
t}\sigma\big(\frac{\eta^2}{4\kappa}(e^{\kappa
t}-1),\frac{x}{s_0}\big)
\end{equation}
where $\sigma^2(t,x)=\mathbb{E}[R_t^2|\exp(I_t-\frac{1}{2}A_t)=x]$
and $R_t$ is a Bessel Process of dimension
$\delta=\frac{4\kappa\theta}{\eta^2}$ starting from $V_0$.
\end{locheston}

\begin{proof}
First, one has the Gy\"{o}ngy volatility formula:
\begin{equation}
\tilde{\sigma}^2(t,x)=\frac{\eta^2}{4}e^{2\kappa
t}\mathbb{E}\big[V_t|S_t=x \big]
\end{equation}
Then using Lemma 2.4, one easily obtains the result.
\end{proof}

\newtheorem{locrem}[green]{Remark}
\begin{locrem}
Let us note that we only have closed-form formulas in cases where
$V_0=0$ and that otherwise we have to go through Laplace transform
inversion techniques.
\end{locrem}

One can propose a general framework for constructing stochastic
volatility models based on Bessel processes. Local volatilities can
be computed through the proposition below whose proof is left to the
reader.

\newtheorem{locbessel2}[green]{Proposition}
\begin{locbessel2}
Let us consider the following stochastic volatility model:
\begin{eqnarray*}
\frac{dS_t}{S_t}&=&\sqrt{v_t}d\beta_t,\indent S_{\{t=0\}}=S_0\\
v_t&=&\frac{g'(t)}{f(t)}V_t,\indent v_0=\frac{g'(0)}{f(0)}V_0 \\
dV_t&=&\bigg(\delta f(t)g'(t)+ \frac{f'(t)}{f(t)} V_t\bigg)dt+\sqrt{f(t)}g'(t)\sqrt{V_t}dW_t\\
d<\beta_{.},W_{.}>_t&=&\rho dt
\end{eqnarray*}
where $\beta_t$ and $W_t$ are Brownian motions, $f$ is a positive
continuously differentiable function and $g$ an increasing
$\mathcal{C}^1$ function.
\\Then the local volatility $\tilde{\sigma}$ that gives us the expected mimicking properties, satisfies the
following equation:
\begin{equation}
\tilde{\sigma}(t,x)=g'(t)\sigma\big(g(t),\frac{x}{S_0}\big)
\end{equation}
where $\sigma^2(t,x)=\mathbb{E}[R_t^2|\exp(I_t-\frac{1}{2}A_t)=x]$
and $R_t$ is a Bessel Process of dimension $\delta$ starting from
$V_0$.
\end{locbessel2}

\section{Pricing Equity Derivatives under Stochastic Interest Rates}
\subsection{A Local Volatility Framework}

With the growth of hybrid products, it has been necessary to take
properly into account the stochasticity of interest rates in FX or
Equity models in a way that makes the equity volatility surface
calibration easy at a given interest rate parametrization. It has
been now a while that people have been considering interest rates as
stochastic for long-dated Equity or FX options, but they have not
been thinking about it in terms of calibration issues. Besides,
according to the interest rates part of an equity - interest rates
hybrid product for example, the instruments on which the interest
rates model will be calibrated are different; hence it becomes
necessary to parameterize the volatility surface efficiently. For
most of hybrid products, no forward volatility dependence is
involved and then a local volatility framework is sufficient. Let us
now consider a local volatility model with stochastic interest
rates:
\begin{equation*}
\frac{dS_t}{S_t}=r_t dt +\sigma(t,S_t) dW_t
\end{equation*}
where $r_t$ is a stochastic process and $\sigma$ a deterministic function.\\
Now, we can observe that equation (\ref{proofvol}) is still valid
under stochastic rates and we may then write
\begin{equation*}
d_t C(t,x)=x\mathbb{E}[r_t e^{-\int_0^t r_s
ds}\mathbf{1}_{\{S_t>x\}}] dt+\frac{1}{2}\lim_{\epsilon\rightarrow
0}\mathbb{E}[\frac{1}{\epsilon}\mathbf{1}_{\{x\leq S_t<
x+\epsilon\}}e^{-\int_0^t r_s ds} \sigma^2(t,S_t) S^2_t] dt
\end{equation*}
The second term of the right-hand side may be written as follows
\begin{equation*}
\mathbb{E}[e^{-\int_0^t r_s ds}\sigma^2(t,S_t)
S^2_t\mathbf{1}_{\{x\leq S_t<
x+\epsilon\}}]=\mathbb{E}\big[\mathbb{E}[e^{-\int_0^t r_s
ds}|S_t]\sigma^2(t,S_t) S_t^2\mathbf{1}_{\{x\leq S_t<
x+\epsilon\}}\big]
\end{equation*}
and then we have
\begin{equation*}
\lim_{\epsilon\rightarrow
0}\frac{1}{\epsilon}\mathbb{E}[e^{-\int_0^t r_s ds}\sigma^2(t,S_t)
S^2_t\mathbf{1}_{\{x\leq S_t< x+\epsilon\}}] =
x^2\sigma^2(t,x)q_t(x)\mathbb{E}[e^{-\int_0^t r_s ds}|S_t=x]
\end{equation*}
where $q_t(x)$ is the value of the density of $S_t$ in $x$. It is
easily shown as well that
\begin{equation*}
\frac{\partial^2 C}{\partial x^2}=q_t(x)\mathbb{E}[e^{-\int_0^t r_s
ds}|S_t=x]
\end{equation*}
Let us now define the $t$-forward measure $\mathbb{Q}^t$ (see Geman
(1989), Jamshidian (1989)) by
\begin{equation*}
\frac{d\mathbb{Q}^t}{d\mathbb{Q}}=\frac{e^{-\int_0^t r_s
ds}}{B(0,t)}\indent \textmd{where}\indent
B(0,t)=\mathbb{E}[e^{-\int_0^t r_s ds}]
\end{equation*}
Hence, we finally obtain an extension of Dupire (1994)'s formula :
\begin{equation*}
\sigma^2(t,x)=\frac{\frac{\partial C}{\partial
t}-xB(0,t)\mathbb{E}^t[r_t\mathbf{1}_{\{S_t>x\}}]}{\frac{x^2}{2}\frac{\partial^2
C}{\partial x^2}}
\end{equation*}
Under a $T$-forward measure for $T\geq t$, one has
\begin{equation*}
\mathbb{E}^T[r_T|\mathcal{F}_t]=f(t,T)
\end{equation*}
where $f(t,T)$ is the instantaneous forward rate. To conclude this
subsection, we can first notice that this slight extension of Dupire
equation may be also written
\begin{equation}
\sigma^2(t,x)=\frac{\frac{\partial C}{\partial
t}+xf(0,t)\frac{\partial C}{\partial x}-xB(0,t)\mathbb{C}
ov^t(r_t;\mathbf{1}_{\{S_t>x\}})}{\frac{x^2}{2}\frac{\partial^2
C}{\partial x^2}}
\end{equation}
We then assume that it is possible to extract from markets prices
the quantities $\mathbb{C}ov^t(r_t;\mathbf{1}_{\{S_t>x\}})$ (i.e.\
there exist tradeable assets from which we could obtain these
covariances) in order to add stochastic interest rates to the usual
local volatility framework. For the remainder of the paper, we
denote this assumption the $(HC)$-Hypothesis that stands for Hybrid
Correlation hypothesis. Under this market hypothesis, one is able to
calibrate a local volatility surface with stochastic interest rates
implied by the derivatives' market prices.

\subsection{Mimicking Stochastic Volatility Models}

In this subsection, we consider the case of a stochastic volatility
model with stochastic interest rates and see how it is possible to
connect it to a local volatility framework. Let us consider the
following diffusion
\begin{equation*}
\frac{dS_t}{S_t}=r_t dt +\sqrt{V_t} dW_t
\end{equation*}
with $V_t$ a stochastic process and let us use equation
(\ref{proofvol}) in order to exhibit a new mimicking property:
\begin{equation*}
d_t C(t,x)=x\mathbb{E}[r_t e^{-\int_0^t r_s
ds}\mathbf{1}_{\{S_t>x\}}] dt+\frac{1}{2}\lim_{\epsilon\rightarrow
0}\mathbb{E}[\frac{1}{\epsilon}\mathbf{1}_{\{x\leq S_t<
x+\epsilon\}}e^{-\int_0^t r_s ds} V_t S^2_t] dt
\end{equation*}
Then,
\begin{eqnarray*}
\lim_{\epsilon\rightarrow
0}\frac{1}{\epsilon}\mathbb{E}[e^{-\int_0^t r_s ds}V_t
S^2_t\mathbf{1}_{\{x\leq S_t< x+\epsilon\}}] &=&
\lim_{\epsilon\rightarrow
0}\frac{1}{\epsilon}\mathbb{E}\big[\mathbb{E}[V_t e^{-\int_0^t r_s
ds}|S_t] S_t^2\mathbf{1}_{\{x\leq S_t< x+\epsilon\}}\big]
\\&=& x^2q_t(x)\mathbb{E}[V_te^{-\int_0^t r_s
ds}|S_t=x]\\
&=&x^2 \frac{\partial^2 C}{\partial
x^2}\frac{\mathbb{E}[V_te^{-\int_0^t r_s
ds}|S_t=x]}{\mathbb{E}[e^{-\int_0^t r_s ds}|S_t=x]}
\end{eqnarray*}
Hence, we obtain
\begin{equation*}
\frac{\mathbb{E}[V_te^{-\int_0^t r_s
ds}|S_t=x]}{\mathbb{E}[e^{-\int_0^t r_s
ds}|S_t=x]}=\frac{\frac{\partial C}{\partial
t}+xf(0,t)\frac{\partial C}{\partial x}-xB(0,t)\mathbb{C}
ov^t(r_t;\mathbf{1}_{\{S_t>x\}})}{\frac{x^2}{2}\frac{\partial^2
C}{\partial x^2}}
\end{equation*}

\paragraph{Spot Mimicking Property}
Finally, if there exists a stochastic process, solution of the
following SDE
\begin{equation*}
\frac{dX_t}{X_t}=r_t dt +\sigma(t,X_t) dW_t
\end{equation*}
such that the one-dimensional marginals of the triple $(r_t,\int_0^t
r_s ds,X_t)$ are the same as $(r_t,\int_0^t r_s ds,S_t)$, then by
identification one must have
\begin{equation*}
\sigma^2(t,x)=\frac{\mathbb{E}[V_te^{-\int_0^t r_s
ds}|S_t=x]}{\mathbb{E}[e^{-\int_0^t r_s ds}|S_t=x]}
\end{equation*}
The existence is easily proven in the cases where $(r_t, t\geq0)$ is
a Markovian diffusion. Hence, we exhibit a strong mimicking property
since we obtained an explicit way to construct a local volatility
surface.

\newtheorem{detrates}[green]{Remark}
\begin{detrates}
We may notice that if interest rates are deterministic, we recover
the well-known formula (\ref{mimic}).
\end{detrates}

\paragraph{Forward Mimicking Property}
Let us now write a Forward mimicking property by applying
Gy\"{o}ngy's result to match the one dimensional marginals of a
stochastic volatility model and of a local volatility one:\\ If one
defines $F^{(1)}_t=S_t e^{-\int_0^t r_s ds}$ and $F^{(2)}_t=X_t
e^{-\int_0^t r_s ds}$ where $S$ and $X$ are defined above, we obtain
the existence of diffusions $Y^{(1)}_t$ and $Y^{(2)}_t$ solutions of
\begin{equation*}
\frac{dY^{(i)}_t}{Y^{(i)}_t}=\Sigma_{(i)}(t,Y^{(i)}_t) dW_t
\end{equation*}
for $i=1,2$ such as
\begin{eqnarray*}
\Sigma_{(1)}^2(t,x)&=&\mathbb{E}[V_t|S_t=x e^{\int_0^t r_s
ds}]\\
\Sigma_{(2)}^2(t,x)&=&\mathbb{E}[\sigma^2(t,x e^{-\int_0^t r_s
ds})|X_t=x e^{\int_0^t r_s ds}]
\end{eqnarray*}
Since the one-dimensional marginals of $F^{(1)}_t$ and $F^{(2)}_t$
must be equal, one obtains
\begin{equation}
\mathbb{E}[V_t|S_t=x e^{\int_0^t r_s ds}]=\mathbb{E}[\sigma^2(t,x
e^{-\int_0^t r_s ds})|X_t=x e^{\int_0^t r_s ds}]
\end{equation}
We consequently obtain an implicit way to construct a local
volatility surface we say that this relation is weak in the sense
that it is a weak mimicking distribution property which is involved
in the above relation.

\subsection{From a Deterministic Interest Rates
Framework to a Stochastic one}

Going from a framework to another is valuable for calibration
issues. Let us assume, for instance that a model has been calibrated
with deterministic interest rates and that one wants to recalibrate
the same model assuming stochastic interest rates. Let us introduce
some notation to define the different kinds of frameworks we will go
through in this subsection.\\

\noindent \textbf{Notation}\\
\textbf{LV} stands for Local Volatility, \textbf{SV} stands for
Stochastic Volatility, \textbf{DIR} stands for Deterministic
Interest Rates and \textbf{SIR} stands for Stochastic Interest Rates

\paragraph{From DIR-LV to SIR-LV}
Let us first consider the local volatility case. Under deterministic
interest rates, the stock price dynamics are driven by the equation
\begin{equation*}
\frac{dS_t}{S_t}=r(t) dt +\sigma(t,S_t)dW_t
\end{equation*}
while under stochastic interest rates it would be
\begin{equation*}
\frac{dS_t}{S_t}=r_t dt +\overline{\sigma}(t,S_t)dW_t
\end{equation*}
and we know that both local volatility functions solve the following
implied equations:
\begin{eqnarray*}
\sigma^2(t,x)&=&\frac{\frac{\partial C}{\partial
t}+xf(0,t)\frac{\partial C}{\partial
x}}{\frac{x^2}{2}\frac{\partial^2 C}{\partial
x^2}}\\
\overline{\sigma}^2(t,x)&=&\frac{\frac{\partial C}{\partial
t}+xf(0,t)\frac{\partial C}{\partial x}-xB(0,t)\mathbb{C}
ov^t(r_t;\mathbf{1}_{\{S_t>x\}})}{\frac{x^2}{2}\frac{\partial^2
C}{\partial x^2}}
\end{eqnarray*}
where $f(0,t)=r(t)$.\\

\indent Now, if the prices involved in the estimation of the local
volatility surfaces are observed on markets and respect the
$(HC)$-Hypothesis, one may write
\begin{equation}
\sigma^2(t,x)-\overline{\sigma}^2(t,x)=\frac{2B(0,t)\mathbb{C}
ov^t(r_t;\mathbf{1}_{\{S_t>x\}})}{x\frac{\partial^2 C}{\partial
x^2}}
\end{equation}
\\
\paragraph{From DIR-SV to SIR-SV}
If we assume that a general It\^{o} process drives the volatility we
will write
\begin{eqnarray*}
\frac{dS^{(1)}_t}{S^{(1)}_t}&=&r(t) dt +\sqrt{V^{(1)}_t}dW_t\\
\frac{dS^{(2)}_t}{S^{(2)}_t}&=&r_t dt +\sqrt{V^{(2)}_t}dW_t
\end{eqnarray*}
and then, if $S^{(1)}$ and $S^{(2)}$ have the same one-dimensional
marginals, we obtain the following relation to relate $V^{(1)}$ to
$V^{(2)}$:
\begin{equation}
\mathbb{E}[V^{(1)}_t|S^{(1)}_t=x]-\frac{\mathbb{E}[V^{(2)}_te^{-\int_0^t
r_s ds}|S^{(2)}_t=x]}{\mathbb{E}[e^{-\int_0^t r_s
ds}|S^{(2)}_t=x]}=\frac{2B(0,t)\mathbb{C}
ov^t(r_t;\mathbf{1}_{\{S^{(2)}_t>x\}})}{x\frac{\partial^2
C}{\partial x^2}}
\end{equation}

\paragraph{From SIR-SV to DIR-LV}
Let us now specify a Heath Jarrow and Morton (1992) diffusion for
the interest rate model and see precisely how one could extract,
using Gy\"{o}ngy's result, the volatility of the forward contract
under deterministic interest rates from the volatility of the
forward contract under stochastic interest rates. Let us recall that
in a standard HJM framework, the instantaneous forward rate follows
\begin{equation*}
df(t,T)=\Big(\sigma(t,T)\int_t^T\sigma(t,u)du\Big)dt+\sigma(t,T)dW^r_t
\end{equation*}
where $\sigma(t,T)$ is a stochastic process adapted to its canonical
filtration and where the price satisfies
\begin{equation*}
B(t,T)=\exp\bigg(-\int_t^T f(t,s) ds\bigg)
\end{equation*}
By definition $r_t=f(t,t)$ and then we obtain
\begin{eqnarray*}
\frac{d B(t,T)}{B(t,T)}&=&r_t dt -\sigma_B(t,T)dW^r_t\\
\sigma_B(t,T)&=&\int_t^T\sigma(t,u)du
\end{eqnarray*}
For our purpose, let us consider a general model
\begin{equation*}
\frac{dS_t}{S_t}=r_t dt +\sqrt{V_t}dW_t
\end{equation*}
then recall the price of the $T$-forward contract written on $S$
\begin{equation*}
F_t^T=\frac{S_t}{B(t,T)}
\end{equation*}
where we assume $d<W,W^r>_t=\rho dt$. We are now able to write the
dynamics of $F_t^T$ under $\mathbb{Q}$ the risk-neutral measure:
\begin{equation*}
\frac{dF_t^T}{F_t^T}=\frac{dS_t}{S_t}-\frac{d<S_{\cdot},B(\cdot,T)>_t}{S_tB(t,T)}-\bigg(\frac{d
B(t,T)}{B(t,T)}-\frac{d<B(\cdot,T)>_t}{B^2(t,T)}\bigg)
\end{equation*}
If we introduce the $T$-forward probability measure as above by
\begin{equation*}
\frac{d\mathbb{Q}^T}{d\mathbb{Q}}=\frac{e^{-\int_0^T r_s
ds}}{B(0,T)}
\end{equation*}
we explain the dynamics of $F_t^T$ under $\mathbb{Q}^T$
\begin{equation*}
\frac{dF_t^T}{F_t^T}=\sqrt{V_t} d\widetilde{W}_t+\sigma_B(t,T)
d\widetilde{W}^r_t
\end{equation*}
where $\widetilde{W}$ and $\widetilde{W}^r$ are Brownian motions
under $\mathbb{Q}^T$ such as
$d<\widetilde{W},\widetilde{W}^r>_t=\rho dt$.\\
We now apply Theorem 2.11 and obtain the existence of a process
$\widetilde{F}_t^T$ solution of an inhomogeneous Markovian
stochastic differential equation
\begin{equation*}
\frac{d\widetilde{F}_t^T}{\widetilde{F}_t^T}=\Sigma_T(t,\widetilde{F}_t^T)d\beta_t
\end{equation*}
where $\beta$ is a Brownian motion and
\begin{equation*}
\Sigma^2_T(t,x)=\mathbb{E}^T[V_t+2\rho\sqrt{V_t}\sigma_B(t,T)+\sigma^2_B(t,T)|S_t=x
B(t,T)]
\end{equation*}
If we consider a local volatility model with deterministic interest
rates as follows
\begin{equation*}
\frac{dS_t}{S_t}=f(0,t) dt +\sigma(t,S_t)d\gamma_t
\end{equation*}
the dynamics of the $T$-forward contract then becomes
\begin{equation*}
\frac{dF_t^T}{F_t^T}=\sigma(t,F_t^T e^{-\int_t^T f(0,s)ds})d\gamma_t
\end{equation*}
and Gy\"{o}ngy's result enables us to conclude that
\begin{equation*}
\Sigma_T(t,x)=\sigma(t,x e^{-\int_t^T f(0,s)ds})
\end{equation*}
Hence, we have proven a new relation that links a local volatility
framework with deterministic interest rates to a stochastic
volatility one with stochastic interest rates, namely
\begin{equation}
\sigma^2(t,x)=\mathbb{E}^T[V_t+2\rho\sqrt{V_t}\sigma_B(t,T)+\sigma^2_B(t,T)|S_t=x
B(t,T) e^{\int_t^T f(0,s)ds}]
\end{equation}
An illustration of this formula can be found for a Black and Scholes
(1973) framework with random rates for example in Hull and White
(1994).

\pagebreak
\section{Conclusion}

This paper recalls well-known results on local volatility and
establishes links to stochastic volatility through the powerful
theorems of Krylov and Gy\"{o}ngy. These general results are then
illustrated with explicit computations of local volatility in
different stochastic volatility models where the volatility process
is a time-space transformation of Bessel processes. In this
framework, we show the impact of the stock-volatility correlation
on the local volatility surface.\\
The local volatility extracted from a stochastic volatility model
allows us to get a precise idea of the skew generated by a
stochastic volatility model. Hence, an important theoretical and
 numerical advantage of generating a local volatility surface from a stochastic
volatility rather than from market option prices is the stability
and the meaningfulness of the surface. Indeed, the local volatility
surface constructed with the Forward PDE equation is known to be
completely unstable whereas as one can see the one built from a
stochastic
volatility is really smooth.\\
With the growth of hybrid products, it has been important to
seriously consider the issue of volatility calibration under
stochastic interest rates and that is the reason why we exhibit
different relations between local volatilities, stochastic
volatilities and derivative prices. It is shown that Dupire (1994)
and Derman and Kani (1998) formulas can easily be extended and that
it is possible to relate any continuous stochastic volatility model
with stochastic interest rates to a local volatility one with
deterministic interest rates. By extending the local volatility
formula to a stochastic rates framework, it is observed that a
market premium for the hybrid correlation risk is to be implied for
the construction of the local volatility surface, which can be
performed under the $(HC)$-Hypothesis as at some point
a market premium for the volatility risk is to be taken into account.\\
A remaining interesting question is the existence of a local
volatility diffusion with a general Ito interest rates process
framework such that the joint law of instantaneous rate, the
discount factor and the stock price is the same as the one in a
stochastic volatility framework.




\end{document}